\documentclass[11pt]{article}
\usepackage{graphicx,subfigure,ifthen,color,algorithm,palatino}
\usepackage{amsmath,amsthm,amssymb,amsfonts,verbatim}
\usepackage{mathrsfs,accents,setspace}
\DeclareMathAlphabet\mathbfcal{OMS}{cmsy}{b}{n}
\newboolean{UnBlinded}
\newboolean{DoubleSpaced}
\setboolean{UnBlinded}{true}
\setboolean{DoubleSpaced}{false}
\def\sumim{\sum_{i=1}^m}
\def\bOmega{\boldsymbol{\Omega}}
\def\bLambda{\boldsymbol{\Lambda}}
\def\bv{\boldsymbol{v}}
\def\bM{\boldsymbol{M}}
\def\bD{\boldsymbol{D}}
\def\vech{\mbox{\rm vech}}
\def\bA{\boldsymbol{A}}
\def\vecof{\mbox{\rm vec}}
\def\bO{\boldsymbol{O}}
\def\AsyVar{\mbox{Asy.Var}}
\def\trigamma{\mbox{trigamma}}
\def\bt{\boldsymbol{t}}
\def\phihat{{\widehat\phi}}
\def\simind{\stackrel{{\tiny \mbox{ind.}}}{\sim}}
\def\dtilde{{\widetilde d}}
\def\bib{\vskip12pt\par\noindent\hangindent=1 true cm\hangafter=1}
\def\bY{\boldsymbol{Y}}
\def\Gsc{{\mathcal G}}
\def\Hsc{{\mathcal H}}
\def\diag{\mbox{diag}}
\def\tr{\mbox{tr}}
\def\bx{\boldsymbol{x}}
\def\argmindum{\mathop{\mbox{\rm argmin}}}
\def\argmin#1{\argmindum_{#1}}
\def\bone{\boldsymbol{1}}
\def\Qsc{{\mathcal Q}}
\def\bI{\boldsymbol{I}}
\def\convprob{\stackrel{P}{\to}}
\def\bK{\boldsymbol{K}}
\def\bL{\boldsymbol{L}}
\def\Psc{{\mathcal P}}
\def\smhalf{{\textstyle{\frac{1}{2}}}}
\def\bU{\boldsymbol{U}}
\def\bzero{\boldsymbol{0}}
\def\bSigma{\boldsymbol{\Sigma}}
\def\bu{\boldsymbol{u}}
\def\argmaxdum{\mathop{\mbox{\rm argmax}}}
\def\argmax#1{\argmaxdum_{#1}}
\def\dA{d_{\mbox{\tiny A}}}
\def\intdA{\int_{\real^{\dA}}}
\def\bQsc{\pmb{\Qsc}}
\def\gothicc{\mathfrak{c}}
\def\gothicg{\mathfrak{g}}
\def\gothich{\mathfrak{h}}
\def\bAsc{\mathbfcal{A}}

\def\HscrdAAAi{\Hsc'_{\mbox{\scriptsize AAA}i}}

\def\dB{d_{\mbox{\tiny B}}}
\def\bX{\boldsymbol{X}}
\def\bbeta{\boldsymbol{\beta}}
\def\betaA{\beta_{\mbox{\scriptsize A}}}

\def\bbetaAMLE{{\widehat\bbeta}_{\mbox{\scriptsize A}}}
\def\bbetaBMLE{{\widehat\bbeta}_{\mbox{\scriptsize B}}}

\def\bSigmaMLE{{\widehat\bSigma}}

\def\vechbSigma{\vech(\bSigma)}
\def\phiMLE{{\widehat\phi}}
\def\bSigmaZero{\bSigma^0}
\def\phiZero{\phi^0}

\def\bXA{\bX_{\mbox{\scriptsize A}}}
\def\bXB{\bX_{\mbox{\scriptsize B}}}

\def\bXAij{\bX_{\mbox{\scriptsize A}ij}}
\def\bXBij{\bX_{\mbox{\scriptsize B}ij}}
\def\bbetaA{\bbeta_{\mbox{\scriptsize A}}}
\def\bbetaB{\bbeta_{\mbox{\scriptsize B}}}

\def\tinybbetaA{\bbeta_{\mbox{\tiny A}}}
\def\tinybbetaB{\bbeta_{\mbox{\tiny B}}}

\def\bbetaAzero{\bbetaA^0}
\def\bbetaBzero{\bbetaB^0}
\def\bSigmaZero{\bSigma^0}

\def\GscrAi{\Gsc_{\mbox{\scriptsize A}i}}

\def\HscrAAi{\Hsc_{\mbox{\scriptsize AA}i}}

\def\linXABUizero{(\bbetaAzero+\bU_i)^T\bXAij+(\bbetaBzero)^T\bXBij}

\def\linXAB{(\bbetaA)^T\bXAij+(\bbetaB)^T\bXBij}
\def\linXABUi{(\bbetaA+\bU_i)^T\bXAij+\bbetaB^T\bXBij}
\def\linXABu{(\bbetaA+\bu)^T\bXAij+\bbetaB^T\bXBij}

\def\linXABU{(\bbetaA+\bU)^T\bXA+\bbetaB^T\bXB}
\def\linXABUzero{(\bbetaAzero+\bU)^T\bXA+(\bbetaBzero)^T\bXB}
\def\mainTheoremNum{1}

\def\sumjni{\sum_{j=1}^{n_i}}

\def\lambdaMin{\lambda_{\mbox{\rm \scriptsize min}}}

\def\dispsumjni{{\displaystyle\sumjni}}

\def\VertF{\Vert_{\mbox{\tiny $F$}}}

\def\BigVertF{\Big\Vert_{\mbox{\tiny $F$}}}

\def\real{{\mathbb R}}
\def\convdist{\stackrel{\mathcal D}{\longrightarrow}}

\def\betazZero{\beta_0^0}

\def\sigsqZero{(\sigma^2)^0}

\def\phiMLE{{\widehat\phi}}
\def\phiZero{\phi^0}

\def\oon{\frac{1}{n}}

\def\myand{\&\ }

\def\betazZero{\beta_0^0}

\def\sigsqZero{(\sigma^2)^0}

\def\pYiGXi{p_{\mbox{\tiny $\bY_i|\bX_i$}}}

\def\dispsumjn{{\displaystyle\sum_{j=1}^{n_i}}}

\def\dAwindow{\dA^{\boxplus}}
\def\psiMLE{{\widehat\psi}}
\def\psiZero{\psi^0}

\def\delVechSigma{\nabla_{\scriptsize\mbox{$\vech(\bSigma)$}}}
\def\delbbeta{\nabla_{\scriptsize\mbox{$\bbeta$}}}
\def\betazZero{\beta_0^0}
\def\bbetaBZero{\bbetaB^0}

\newtheorem{theorem}{\textbf{Theorem}}
\begin{document}

\ifthenelse{\boolean{DoubleSpaced}}{\setstretch{1.5}}{}

\vskip5mm
\centerline{\Large\bf Dispersion Parameter Extension of Precise}
\vskip1mm
\centerline{\Large\bf Generalized Linear Mixed Model Asymptotics}
\vskip5mm
\centerline{\normalsize\sc Aishwarya Bhaskaran and Matt P. Wand}
\vskip5mm
\centerline{\textit{University of Technology Sydney}}
\vskip6mm
\centerline{10th August, 2022}

\vskip6mm

\centerline{\large\bf Abstract}
\vskip2mm

We extend a recently established asymptotic normality theorem for generalized linear mixed
models to include the dispersion parameter. The new results show that the maximum likelihood 
estimators of all model parameters have asymptotically normal distributions with asymptotic
mutual independence between fixed effects, random effects covariance and dispersion parameters.
The dispersion parameter maximum likelihood estimator has a particularly simple asymptotic 
distribution which enables straightforward valid likelihood-based inference.


\section{Introduction}\label{sec:intro}

We extend the main theorem of Jiang \textit{et al.} (2022) to include conditional
maximum likelihood estimation of the dispersion parameter. The essence of the new findings is that 
the dispersion parameter maximum likelihood estimator has a simple asymptotic normal distribution,
identical to that for the generalized linear model case, which is amenable to practical inference.
Moreover, we establish asymptotic orthogonality between the dispersion parameter and the other
model parameters.
  
\begin{figure}[!h]
\centering
\includegraphics[width=\textwidth]{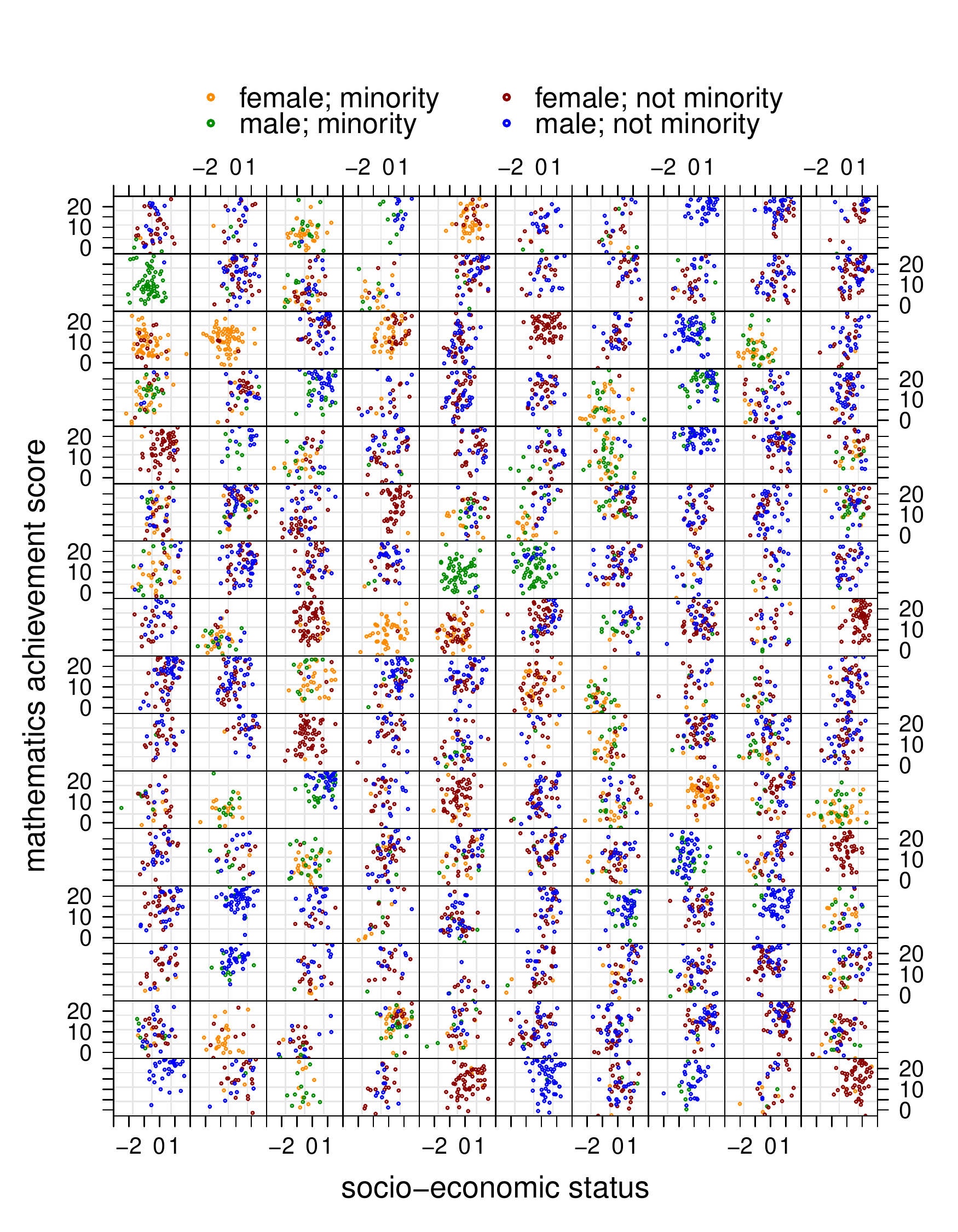}
\caption{\textit{A visualization of the data set concerning mathematics achievement described
in the text. Each panel corresponds to one of $160$ schools and plots a mathematics achievement
score against a socio-economic status index for students in that school. The color-coding corresponds
to the sex and minority status categorical variables.}}
\label{fig:MathAchieveDataCol}
\end{figure}

The motivations and benefits of precise asymptotics for generalized linear mixed models are described
in Section 1 of Jiang \textit{et al.} (2022). As mentioned there, books such as Faraway (2016), 
Jiang \myand Nguyen (2021), McCulloch \textit{et al.} (2008) and Stroup (2013) 
provide summaries and access to the
large literature on generalized linear mixed models.
Figure \ref{fig:MathAchieveDataCol} provides visualization of a data set that potentially benefits
from generalized linear mixed model analysis.  A variable of interest is mathematics 
achievement score for $7,185$ United States of America school 
students across $160$ schools. In Figure \ref{fig:MathAchieveDataCol} each panel corresponds to a 
school and mathematics achievement is plotted against socio-economic status. The color-coding is
according to sex and minority categorical variables. The data are from a 1980s survey known
as ``High School and Beyond'' and an early source is Coleman \textit{et al.} (1982).
The Figure \ref{fig:MathAchieveDataCol} data is stored in the data frame \texttt{MathAchieve} 
within the package \textsf{nlme} (Pinheiro \textit{et al.}, 2022) of the 
\textsf{R} computing environment (\textsf{R} Core Team, 2022).

Estimation of the dispersion parameter, denoted here by $\phi$, was not considered by 
Jiang \textit{et al.} (2022) for a few reasons. One is that $\phi$ is often treated
as a nuisance parameter, with the fixed and random effects being of primary interest.
Another is that maximum likelihood does not apply for the most common families: 
Bernoulli and Poisson. Instead, for these two families, quasi-likelihood is required 
for the $\phi\ne1$ extension. A third reason is that, even when 
maximum likelihood estimation is available for situations such as Gamma response models, it 
is common to use the simpler Pearson estimator instead. Nevertheless, maximum likelihood 
estimation of $\phi$ is viable in important generalized linear mixed model contexts and a 
full treatment of precise asymptotics requires its inclusion.

Section \ref{sec:modelMLE} describes the data, generalized linear mixed models set-up and maximum
likelihood estimation of the model parameters. The focal point of this article is the asymptotic
normality theorem stated in Section \ref{sec:mainResult} with proof provided in supplemental material.
Asymptotically valid likelihood-based inference for the dispersion parameter is discussed in 
Section \ref{sec:validInference}. Sections \ref{sec:modelMLE}--\ref{sec:validInference} are 
confined to reproductive exponential families, which covers the common cases arising in 
generalized linear mixed model applications. The extension to general exponential families
is discussed in Section \ref{sec:genFam}.

\section{Model Description and Maximum Likelihood Estimation}\label{sec:modelMLE}

Consider the class of two-parameter reproductive linear exponential family density functions 
with generic form
\begin{equation}
p(y;\eta,\phi)= \exp[\left\{y\eta - b(\eta)+c(y)\right\}/\phi - d(\phi)-e(y)]h(y)
\label{eq:expDens}
\end{equation}
where $\eta$ is the \emph{natural parameter} and $\phi>0$ is the \emph{dispersion parameter}.
The definition of a reproductive exponential family is given in, for example, J{\o}rgensen (1987).
Not all two-parameter linear exponential family density functions are reproductive. However, those 
families commonly used in applications of generalized linear and mixed models are reproductive 
and Sections \ref{sec:modelMLE}--\ref{sec:validInference} are restricted to this case. 
We discuss the general case briefly in Section \ref{sec:genFam}.

Table \ref{tab:bcdeh} gives some explicit examples of the $b$, $c$, $d$, $e$ and $h$ functions 
appearing in (\ref{eq:expDens}). In Table \ref{tab:bcdeh} the following notation is used: 
$I(\Psc)=1$  if the condition $\Psc$ is  true and $I(\Psc)=0$ if $\Psc$ is false. 
Theoretical results given in Bl{\ae}sild \myand Jensen (1985) imply that the three families 
listed in Table \ref{tab:bcdeh} are the \emph{only} possibilities for $p(y;\eta,\phi)$, even 
if $-d(\phi)-e(y)$ is relaxed to be a general bivariate function of $(\phi,y)$. Therefore, 
without loss of generality, we can assume that $p(y;\eta,\phi)$ is one of the three forms 
given by Table \ref{tab:bcdeh}.

\begin{table}[!h]
\begin{center}
\begin{tabular}{lccccc}
\hline\\[-1.8ex]
family & $b(\eta)$ & $c(y)$ &$d(\phi)$ & $e(y)$ & $h(y)$ \\[0.1ex]
\hline\\[-1.5ex]
Gaussian  &  $\smhalf\eta^2$  &  $-\smhalf y^2$ & $\smhalf\log(\phi)$ & $\smhalf\log(2\pi)$  & $1$ \\[2ex]
Gamma     &  $-\log(-\eta)$  &  $\log(y)$ & $\log\big(\phi^{1/\phi}\,\Gamma(1/\phi)\big)$   
& $\log(y)$ & $I(y>0)$ \\[2ex]   
Inverse Gaussian &  $-(-2\eta)^{1/2}$  &  $-1/(2y)$ & $\smhalf\log(\phi)$
& $\smhalf\log(2\pi y^3)$ & $I(y>0)$ \\[1ex]   
\hline
\end{tabular}
\end{center}
\caption{\textit{Specific two-parameter exponential families and their
 $b$, $c$, $d$, $e$ and $h$ functions.}}
\label{tab:bcdeh}
\end{table}

In this article we study generalized linear mixed models of the form, for 
observations of the random triples 
$(\bXAij,\bXBij,Y_{ij})$, $1\le i\le m$, $1\le j\le n_i$,
\begin{equation}
\begin{array}{l}
Y_{ij}|\bXAij,\bXBij,\bU_i\ \mbox{independent having density function (\ref{eq:expDens}) with natural} \\[1ex]
\mbox{parameter}\ \linXABUizero\ \mbox{such that the}\ \bU_i\ 
\mbox{are independent}\\[1ex]
\ N(\bzero,\bSigmaZero)\ \mbox{random vectors}.
\end{array}
\label{eq:theModel}
\end{equation}
The $\bU_i$ are $\dA\times1$ unobserved random effects vectors. 
The $\bXAij$ are $\dA\times1$ random vectors corresponding to predictors
that are partnered by a random effect. The $\bXBij$ are $\dB\times1$ random vectors
corresponding to predictors that have a fixed effect only. Let $\bX_{ij}\equiv(\bXAij^T,\bXBij^T)^T$
denote the combined predictor vectors. We assume that the $\bX_{ij}$ and $\bU_i$, for $1\le i\le m$ 
and $1\le j\le n_i$, are totally independent, with the $\bX_{ij}$ each having the
same distribution as the $(\dA+\dB)\times 1$ random vector $\bX=(\bXA^T,\bXB^T)^T$ 
and the $\bU_i$ each having the same distribution as the random vector $\bU$.

For any $\bbetaA$ $(\dA\times1)$, $\bbetaB$ $(\dB\times1)$, $\bSigma$ $(\dA\times\dA)$ that is 
symmetric and positive definite and $\phi>0$, conditional on the $\bX_{ij}$ data, the log-likelihood is 
\begin{equation}
{\setlength\arraycolsep{0pt}
\begin{array}{rcl}
&&\ell(\bbetaA,\bbetaB,\bSigma,\phi)=-\frac{m}{2}\log|2\pi\bSigma|\\[1ex]
&&\ \ +{\displaystyle\sum_{i=1}^m\sum_{j=1}^{n_i}}
\left[\{Y_{ij}(\bbetaA^T\bXAij+\bbetaB^T\bXBij)+c(Y_{ij})\}/\phi
-d(\phi)-e(Y_{ij})\right]\\
&&\ \ +{\displaystyle\sum_{i=1}^m}\log{\displaystyle\int_{\real^{\dA}}}\exp\Bigg[{\displaystyle\sum_{j=1}^{n_i}}
\{Y_{ij}\bu^T\bXAij-b\big(\linXABu\big)\}/\phi-\smhalf\bu^T\bSigma^{-1}\bu\Bigg]\,d\bu.
\end{array}
}
\label{eq:MrsMurray}
\end{equation}
The conditional maximum likelihood estimator of $(\bbetaAzero,\bbetaBzero,\bSigmaZero,\phiZero)$ is 
$$(\bbetaAMLE,\bbetaBMLE,\bSigmaMLE,\phiMLE)=\argmax{\tinybbetaA,\tinybbetaB,\bSigma,\phi}
\ell(\bbetaA,\bbetaB,\bSigma,\phi).$$

Even though conditional maximum likelihood provides a natural estimator for $\phi$ based on
data modeled according to (\ref{eq:theModel}), the relevant literature and software is
such that alternative approaches are common. For the generalized linear model special case,
Section 8.3.6 of McCullagh \myand Nelder (1989) expresses a preference for the natural 
moment-based estimator of $\phi$, which is sometimes referred to as the Pearson estimator.
The functions \textsf{glm()} for generalized linear models and \textsf{glmer()} within the 
\textsf{R} package \textsf{lme4} (Bates \textit{et al.}, 2015) for generalized linear
mixed models each use Pearson estimation of the dispersion parameter. In Section 2 of 
Cordeiro \myand McCullagh (1991) some alternative estimators for the dispersion parameters
are proposed, motivated by bias correction and computational convenience considerations. 
Jo \myand Lee (2017) compare the efficiencies of dispersion parameter estimators for
Gamma generalized linear models and recommend conditional maximum likelihood
estimation compared with the Pearson and Cordeiro \myand McCullagh (1991) estimators.

\section{Main Result}\label{sec:mainResult}

Given the addition of $\phi$ to the set of parameters being estimated, compared with
the Jiang \textit{et al.} (2022) set-up, our aim is to extend the asymptotic normality
result for this enlarged estimation problem. We start by repeating definitions and
conditions from Section 3 Jiang of \textit{et al.} (2022). Let
$$n\equiv{\displaystyle\frac{1}{m}\sumim} n_i=\mbox{average of the within-group sample sizes},$$
$$\bOmega_{\bbetaB}(\bU)\equiv 
E\left\{
b''\Big(\linXABUzero\Big)
\left[\begin{array}{cc}    
\bXA\bXA^T & \bXA\bXB^T \\[1ex]
\bXB\bXA^T & \bXB\bXB^T
\end{array}
\right]
\Bigg|\bU
\right\},
$$
$$\bLambda_{\bbetaB}\equiv
\Bigg(E\Big[\big\{\mbox{lower right $\dB\times\dB$ block of $\bOmega_{\bbetaB}(\bU)^{-1}$}\big\}^{-1}\Big]\Bigg)^{-1}
$$
and $\Vert\bv\Vert\equiv(\bv^T\bv)^{1/2}$ denote the Euclidean norm of a column vector $\bv$.
For a symmetric matrix $\bM$ let $\lambdaMin(\bM)$ denote the smallest eigenvalue of $\bM$.
Also, let $\bD_d$ denote the matrix of zeroes and ones such that $\bD_d\vech(\bA)=\vecof(\bA)$
for all $d\times d$ symmetric matrices $\bA$. The Moore-Penrose inverse of $\bD_d$ is
$\bD_d^+=(\bD_d^T\bD_d)^{-1}\bD_d^T$. Let $d'$ and $d''$ denote the first and second derivatives 
of the $d$ function.

The theorem relies on the following assumptions:

\begin{itemize}
\item[]
\begin{itemize}
\item[(A1)] The number of groups $m$ diverges to $\infty$.
\item[(A2)] The within-group sample sizes $n_i$ diverge to $\infty$ in such a way that
$n_i/n\to C_i$ for constants $0<C_i<\infty$, $1\le i\le m$. Also, $n/m\to0$ as $m$ and $n$ diverge.
\item[(A3)] The distribution of $\bX$ is such 
that 
\begin{equation}
\displaystyle{E\left[\frac{
E\Big[\max\big(1,\Vert\bX\Vert\big)^8\,\max\big\{1,b''\big(\linXABU\big)\big\}^4\Big|\bU\Big]}
{\min\big\{1,\lambdaMin\big(E\{\bXA\bXA^T
\,b''\big(\linXABU\big)|\bU\}\big)\big\}^2}\right]}<\infty
\label{eq:StoneTemplePilots}
\end{equation}
for all $\bbetaA\in\real^{\dA}$, $\bbetaB\in\real^{\dB}$ and $\bSigma$ a $\dA\times\dA$
symmetric and positive definite matrix.
\end{itemize}
\end{itemize}

\begin{theorem}
Assume that conditions (A1)--(A3) hold. Then
%
%
$$
\sqrt{m}\left[
{\setlength\arraycolsep{0pt}
\begin{array}{c}
\bbetaAMLE-\bbetaAzero\\[1.5ex]
\sqrt{n}\Big(\bbetaBMLE-\bbetaBzero\Big)\\[1.5ex]
\vech(\bSigmaMLE-\bSigmaZero)\\[1.5ex]
\sqrt{n}(\phiMLE-\phiZero)
\end{array}
}
\right]\convdist N\left(
\left[
{\setlength\arraycolsep{0pt}
\begin{array}{c}
\bzero\\[2.5ex]
\bzero\\[2.5ex]
\bzero\\[2.5ex]
0
\end{array}
}
\right],
\left[
{\setlength\arraycolsep{0pt}
\begin{array}{cccc}
\bSigmaZero & \bO                    &    \bO  & \bO \\[1ex]
\bO         & \phiZero\bLambda_{\bbetaB} &    \bO  & \bO\\[1ex]
\bO         & \bO                    & 2\bD_{\dA}^{+}(\bSigmaZero\otimes\bSigmaZero)\bD_{\dA}^{+T} & \bO\\[2ex]
\bO         & \bO                    & \bO  & \begin{Large}\displaystyle{\frac{1}{2d'(\phiZero)/\phiZero+d''(\phiZero)}}
\end{Large}
\end{array}
}
\right]
\right).
$$
\label{thm:mainTheorem}
\end{theorem}
\vskip1mm\noindent
A proof of Theorem \mainTheoremNum\ is given in Section \ref{sec:theProof} of
the supplemental material.

Some remarks concerning Theorem \mainTheoremNum\ are:
\begin{itemize}
\item[1.] The diagonal blocks appearing in the Multivariate Normal covariance matrix
of Theorem 1 correspond to asymptotic covariances of: 
(1) fixed effects partnered by a random effect, (2) fixed effects not partnered by a 
random effect, (3) random effects covariance parameters and (4) the dispersion 
parameter. For the first and third of these types of parameters
the asymptotic variances have order $m^{-1}$. The remaining parameters,
including the dispersion parameter, have order $(mn)^{-1}$ asymptotic
variances.
\item[2.] Theorem 1 reveals asymptotic orthogonality between 
$\phi$ and $(\bbetaA,\bbetaB,\bSigma)$. This is in addition to the asymptotic
orthogonality between the components of $(\bbetaA,\bbetaB,\bSigma)$,
established by Theorem 1 Jiang \textit{et al.} (2022).
\item[3.] The asymptotic distribution of $\phiMLE$ is the same as that 
arising for the generalized linear model special case of (\ref{eq:theModel})
when there are no random effects. In other words, the asymptotic behaviour
of $\phiMLE$ is not impacted by the extension from generalized linear models
to generalized linear mixed models.
\item[4.] After obtaining the $2d'(\phi)/\phi+d''(\phi)$ expressions for the 
specific $d$ functions of Table \ref{tab:bcdeh} and simplifying, the asymptotic 
variances of $\phiMLE$ become:
\begin{equation}
\AsyVar(\phiMLE)
=\left\{
\begin{array}{cl}
\displaystyle{\frac{2(\phiZero)^2}{mn}} & \mbox{for the Gaussian and the}\\[0ex]
                                        & \mbox{Inverse Gaussian families},\\[1ex]
\displaystyle{\frac{(\phiZero)^4}{\{\mbox{trigamma}(1/\phiZero)-\phiZero\}\,mn}}
& \mbox{for the Gamma family.} 
\end{array}
\right.
\label{eq:KegBrewPub}
\end{equation}

\item[5.] The $\mbox{trigamma}$ function has the following asymptotic expansion:
$$\trigamma(x)=\frac{1}{x}+\frac{1}{2x^2}+\frac{1}{6x^3}-\frac{1}{30 x^5}+\ldots.$$
where the coefficients are simple functions of Bernoulli numbers. It follows that
$$\frac{(\phiZero)^4}{\mbox{trigamma}(1/\phiZero)-\phiZero}\approx 2(\phiZero)^2
\quad\mbox{for small values of}\ \phiZero.
$$
This connects the asymptotic variance results of (\ref{eq:KegBrewPub}) for low 
values of the dispersion parameter.
\item[6.] The moment condition in (A3) is sufficient but not necessary for Theorem 1
to hold. A recently discovered erratum has resulted in the replacement of the second
moment in the second numerator factor in (\ref{eq:StoneTemplePilots}) by the fourth moment,
compared with (A3) of Jiang \textit{et al.} (2022). For various special cases of 
(\ref{eq:expDens}) and (\ref{eq:theModel}) mathematical analysis can be used to 
replace (A3) by a simpler moment condition. The Gaussian response case is such
that $b''=1$ and relatively simple arguments show that (A3) can be replaced by 
$$E(\Vert\bX\Vert^8)<\infty\quad\mbox{and}\quad\mbox{no entry of $\bXA$ is the 
zero degenerate random variable}.$$
In this case it is clear that (A3) holds for sufficiently light-tailed $\bX$ distributions
but fails if, for example, an entry of $\bXB$ has a heavy-tailed distribution such as the 
$t$ distribution with a low degrees of freedom parameter.
A more involved example involves Poisson responses with $\bXA=1$, corresponding to random
intercepts. In Section \ref{sec:IllusSimplic} of the supplemental material it 
is shown that (A3) can be replaced by the moment generating function existence condition
$$E\{\exp(\bt^T\bXB)\}<\infty\ \mbox{for all $\bt\in\real^{\dB}$}$$
which is satisfied by, for example, $\bXB$ having a Multivariate Skew-Normal
distribution (Azzalini \myand Dalla Valle, 1996). The arguments in Section \ref{sec:IllusSimplic}
of the supplemental material give a flavor of what is involved to simplify (A3). Other cases 
require more detailed mathematical analysis.
\end{itemize}

\section{Asymptotically Valid Inference}\label{sec:validInference}

An immediate consequence of Theorem 1 is
%
%
$$
\sqrt{mn}\,
\big\{2d'(\phiZero)/\phiZero+d''(\phiZero)\big\}^{1/2}\Big(\phiMLE-\phiZero\Big)
\convdist N(0,1).
$$
This asymptotic normality result still holds when the unknown
quantities on the left-hand side are replaced by consistent estimators,
often referred to as \emph{studentization}. Hence an asymptotically valid $100(1-\alpha)\%$ 
confidence interval for $\phiZero$ is 
%
%
\begin{equation}
\phiMLE\pm\Phi^{-1}(1-\smhalf\alpha)
\Big[\big\{2d'(\phiMLE)/\phiMLE + d''(\phiMLE)\big\}mn\Big]^{-1/2}
\label{eq:generalCI}
\end{equation}
where $\Phi$ is the $N(0,1)$ cumulative distribution function. It follows that Theorem 1 provides simple 
closed form asymptotically valid inference for $\phiZero$. For the specific families listed 
in Table \ref{tab:bcdeh}, (\ref{eq:generalCI}) simplifies to
$$\phiMLE\pm\Phi^{-1}(1-\smhalf\alpha)\sqrt{\frac{2\,\phihat^{\,2}}{mn}}
\quad\mbox{for the Gaussian and Inverse Gaussian families}
$$
and
\begin{equation}
\phiMLE\pm\Phi^{-1}(1-\smhalf\alpha)\sqrt{\frac{\,\phihat^{\,4}}
{\{\mbox{trigamma}(1/\phiMLE)-\phiMLE\}mn}}
\quad\mbox{for the Gamma family}.
\label{eq:CrisisWhatCrisis}
\end{equation}
%

\begin{figure}[!h]
\centering
\includegraphics[width=0.9\textwidth]{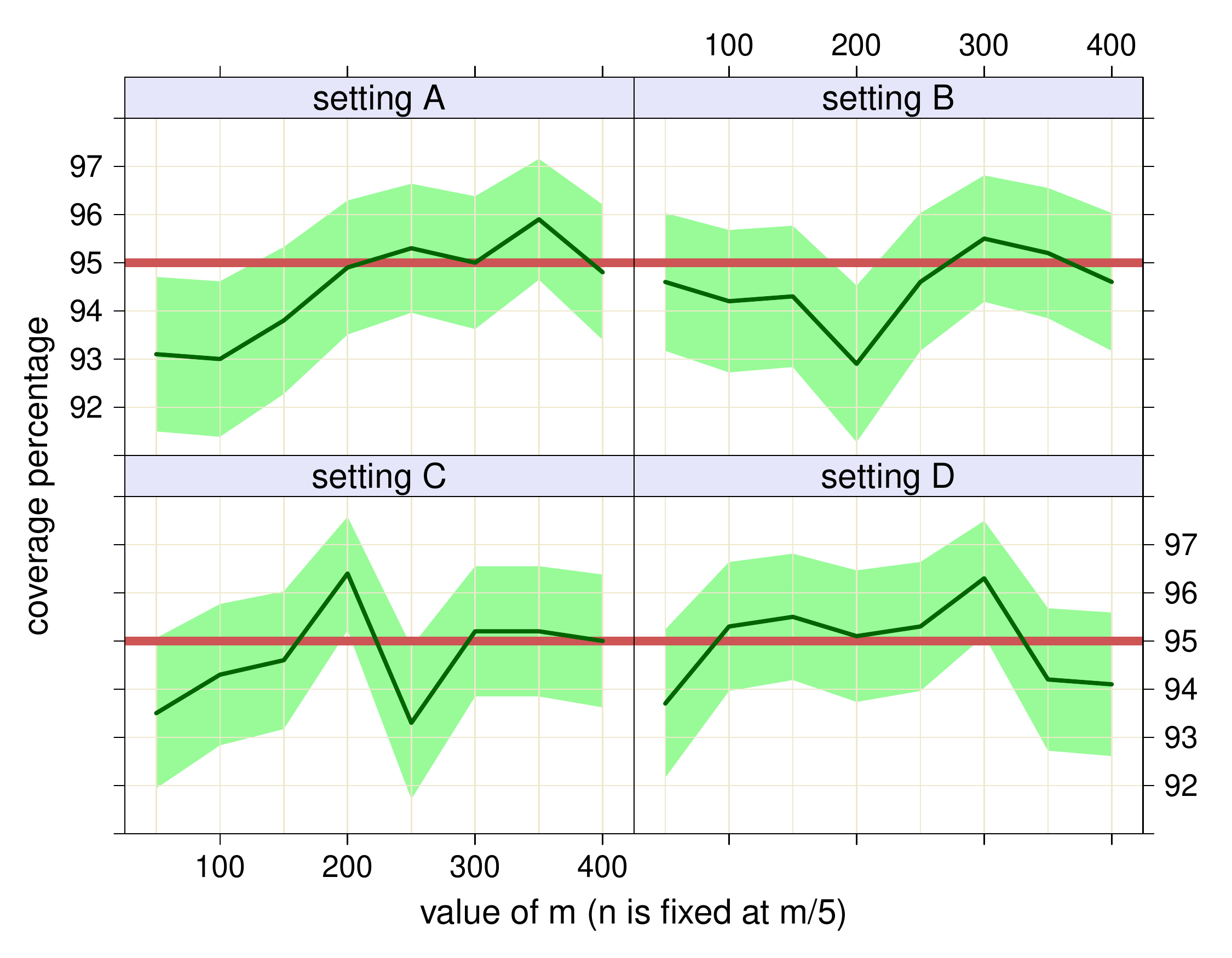}
\caption{\textit{Actual coverage percentage of nominally 95\% confidence
intervals for $\phi^0$ in $\dA=1$ and $\dB=3$ Gamma mixed model with four different
true parameter vectors. The nominal percentage is shown as a thick horizontal line. 
The percentages are based on $1,000$ replications. The values of $m$ are 
$50,100,150,200,250,300,350,400$. The value of $n$ is fixed at $m/5$. The 
shaded regions correspond to plus and minus two standard
errors of the sample proportions.}}
\label{fig:CIcovPlotPhi}
\end{figure}

A simulation study was run to assess the actual coverages of Theorem 1-based 
confidence intervals for the dispersion parameter for the $\dA=1$ and $\dB=3$ 
Gamma mixed model. In this scalar random effects case, $\bbetaA$ and $\bSigma$ are 
replaced by the scalar parameter symbols $\beta_0$ and $\sigma^2$. 
The true parameter vector $(\betazZero,\bbetaBZero,\sigsqZero,\phiZero)$ 
had the following four settings:
\begin{center}
\begin{tabular}{ll}
   setting A: & $(-2.78,-1.55,0,0.98,0.25,0.54)$, \\[0ex]
   setting B: & $(-4.06,-2.41,0.16,-3.93,0.52,1.92)$,\\[0ex]
   setting C: & $(-8.55,3.13,-7.82,-0.23,1.27,0.86)$,\\[0ex]  
   setting D: & $(-14.45,8.78,0.41,-3.32,1.88,2.11)$
\end{tabular}
\end{center}
and the $\bXBij$s were generated independently from the uniform distribution over the unit
cube. The number of groups $m$ varied over the set $\{50,100,150,200,250,300,350,400\}$ 
and the sample size within each group was $n$ fixed at $m/5$.
For each of the possible combinations of the true parameter vector
and the sample size pair we simulated $1,000$ replications. Established generalized
linear mixed model software, such the \texttt{glmer()} function within the 
\textsf{R} package \textsf{lme4} (Bates \textit{et al.}, 2015), does not use maximum
likelihood to estimate $\phi$. For this $\dA=1$ case the likelihood surface (\ref{eq:MrsMurray}) 
was relatively straightforward to evaluate exactly using univariate quadrature 
and maximize using a derivative-free optimization algorithm such as the Nelder-Mead simplex 
method (Nelder \myand Mead, 1965) with \texttt{glmer()} starting values. After obtaining the 
maximum likelihood estimates, we computed 95\% confidence intervals based on 
(\ref{eq:CrisisWhatCrisis}) with $\alpha=0.05$. 

Figure \ref{fig:CIcovPlotPhi} summarizes the empirical coverage results. The shaded
regions around the line segments in Figure \ref{fig:CIcovPlotPhi} indicate plus
and minus two standard errors for the sample proportions.
It is seen that the actual coverage percentages are in keeping with the advertized level
of 95\% across all settings and sample sizes.
There is a tendency for the empirical coverages to get closer, on average, to 95\% 
as $m$ and $n$ in increase, as expected for statistical inference
based on leading term asymptotics.

\section{Analysis of the Mathematics Achievement Data}

As an illustration of the asymptotically valid inference results given in the 
previous section we conducted an analysis of the mathematics 
achievement data shown in Figure \ref{fig:MathAchieveDataCol} using a Gaussian 
response version of (\ref{eq:theModel}). Specifically, we considered the model 
\begin{equation}
{\setlength\arraycolsep{1pt}
\begin{array}{c}
\texttt{MathAchieve}_{ij}|U_{0i},U_{1i}\simind
N\Big(\beta^0_0+U_{0i}+(\beta^0_1+U_{1i})\,\texttt{SES}_{ij}+\beta^0_2\,\texttt{isMale}_{ij}\\
\qquad\qquad\qquad\qquad+\beta_3^0\,\texttt{isMinority}_{ij},\phiZero\Big),\\[1.5ex]
\left[
\begin{array}{c}
U_{0i}\\[1ex]
U_{1i}
\end{array}
\right]\simind N\left(\left[\begin{array}{c} 
0\\[1ex]
0
\end{array}
\right],
\left[
\begin{array}{cc}
(\sigma_0^0)^2           & \rho\sigma_0^0\sigma_1^0\\[1ex]
\rho\sigma_0^0\sigma_1^0 & (\sigma_1^0)^2
\end{array}
\right]
\right),\quad 1\le i\le m,\quad 1\le j\le n_i,\
\end{array}
}
\label{eq:MathAchModel}
\end{equation}
where $\simind$ stands for ``independently distributed as''
and $n_i$ is the number of students in the $i$th school.
In (\ref{eq:MathAchModel}) $\texttt{MathAchieve}_{ij}$ and $\texttt{SES}_{ij}$ denote, 
respectively, the mathematics achievement and socio-economic status scores for the $j$th 
student in the $i$th school. In addition 
$$
\texttt{isMale}_{ij}=\left\{
\begin{array}{rl}
0 &\mbox{if the $j$th student in the $i$th school is female,}\\[0.5ex]
1 &\mbox{if the $j$th student in the $i$th school is male}\\
\end{array}
\right.
$$
and $\texttt{isMinority}_{ij}$ is defined similarly for minority status.
The $n_i$ values range over $14$ to $67$. 

Table \ref{tab:MathAchRes} shows the estimates of the fixed effects
and standard deviation parameters based on conditional maximum likelihood estimation
as described in Section \ref{sec:modelMLE}. 
Also shown in Table \ref{tab:MathAchRes} are approximate 95\% confidence
intervals based on Theorem \ref{thm:mainTheorem} and studentization.
The confidence interval for $\beta_2^0$ indicates a statistically significant 
elevation in mean mathematics achievement score for the male student population.
The confidence interval for $\beta_3^0$ indicates a statistically significant lowering
for minority students. The confidence intervals for $\sigma_0^0$ and $\sigma_1^0$ 
are both within the positive half-line, which indicates significant heterogeneities
in the intercepts and slopes of the social economic status effects.
The last row of Table \ref{tab:MathAchRes} provides an estimate and 95\% confidence
interval for the within-school error standard deviation.
%
\begin{table}[!h]
\begin{center}
\begin{tabular}{crr}
\hline
parameter      & estimate  &   95\% confid. interv.\\
\hline\\[-1.9ex]
$\beta_0^0$    &  12.93    &   $(12.55,13.31)$  \\
$\beta_1^0$    &  2.097    &   $(1.875,2.319)$  \\
$\beta_2^0$    &  1.219    &   $(0.9003,1.537)$  \\
$\beta_3^0$    & $-2.999$  &   $(-3.404,-2.594)$ \\ 
$\sigma_0^0$   &  1.903    &   $(1.682,2.101)$  \\  
$\sigma_1^0$   &  0.4964   &   $(0.4386,0.5481)$ \\
$\sqrt{\phi^0}$&  5.982    &   $(5.883,6.079)$ \\ 
\hline   
\end{tabular}
\end{center}
\caption{\textit{Maximum likelihood estimates and approximate 95\% confidence intervals 
of the fixed effects parameters and standard deviation parameters for the fit of the 
model (\ref{eq:MathAchModel}) to the mathematics achievement data.}}
\label{tab:MathAchRes}
\end{table}

An examination of the residuals revealed reasonable accordance with model assumptions
but some heteroscedasticity. More delicate modelling, involving variance
function extensions of (\ref{eq:theModel}) and other response families may
lead to model fit improvements.

\section{General Two-Parameter Exponential Family Extension}\label{sec:genFam}

We now turn attention to general two-parameter exponential
families. If the restriction to reproductive exponential families is removed then (\ref{eq:expDens}) should
be replaced by 
\begin{equation}
p(y;\eta,\phi)= \exp\Big[\left\{y\eta - b(\eta)+c(y)\right\}/\phi - \dtilde(y,\phi)\Big]h(y)
\label{eq:PuffNStuff}
\end{equation}
where $\dtilde(y,\phi)$ is some bivariate function of $y$ and $\phi$ that is not necessarily
additive in its arguments. Section 2.1 of J{\o}rgensen (1987) describes a procedure for
generating versions of (\ref{eq:PuffNStuff}) from any distribution possessing a moment
generating function. An example of a non-reproductive version of (\ref{eq:PuffNStuff})
given there is such that
$$b(\eta)=\log\Big(-\eta-\sqrt{\eta^2-1}\Big),\ \ c(y)=0,\ \ 
\dtilde(y,\phi)=-\log\Big\{I_{1/\phi}(y/\phi)\big/(y\phi)\Big\}\ \ \mbox{and}\ \ h(y)=I(y>0)$$
where $I_{\nu}$ denotes the modified Bessel function of the first kind with index $\nu$
(e.g. Section 8.431 of Gradshteyn \myand Ryzhik, 1994).

For this more general set-up, the Theorem 1 arguments still apply with $\dtilde(y,\phi)$ replacing
$d(\phi)+e(y)$ and the lower right block of the covariance matrix on the right-hand side
of the Theorem 1 statement generalizes to 
\begin{equation}
(\phiZero)^4\Bigg/
E\left(\left[\frac{\partial^2}{\partial\psi^2} 
\dtilde\left(Y,\frac{1}{\psi}\right)\right]_{\psi=\frac{1}{\phiZero}}
\right)
\label{eq:JebbiePhillips}
\end{equation}
where $Y$ is a random variable having the same distribution as the $Y_{ij}$s. 
Method of moments estimation of the expectation over $Y$ leads to the following
extension of (\ref{eq:generalCI}) for approximate $100(1-\alpha)$\% confidence
intervals for $\phiZero$: 
\begin{equation}
\phiMLE\pm\Phi^{-1}(1-\smhalf\alpha)\phihat^{\,2}
\left(\sumim\sum_{j=1}^{n_i}
\left[\frac{\partial^2 \dtilde\left(Y_{ij},1/\psi\right)}
{\partial\psi^2}\right]_{\small{\mbox{$\psi=1\big/\phiMLE$}}}\right)^{-1/2}.
\label{eq:MaggieTaberer}
\end{equation}
The practical relevance of (\ref{eq:JebbiePhillips}) and (\ref{eq:MaggieTaberer})
is much lower than for the reproductive exponential family special case and
they have been included for completeness.

\section*{Acknowledgements}

We are grateful to Gordon Smyth for advice related to this research. 
This research was partially supported by Australian Research Council
grant DP180100597.

\section*{References}

\bib
Azzalini, A. \myand Dalla Valle, A. (1996). The multivariate
skew-normal distribution. \textit{Biometrika}, \textbf{83}, 715--726.

\bib
Bates, D., Maechler, M., Bolker, B. \myand Walker, S. (2015).
Fitting linear mixed-effects models using \textsf{lme4}. 
\textit{Journal of Statistical Software}, \textbf{67(1)}, 1--48.

\bib
Bl{\ae}sild, P. \myand J.L. Jensen (1985).
Saddlepoint formulas for reproductive exponential models.
\textit{Scandinavian Journal of Statistics}, \textbf{12}, 193--202.

\bib
Coleman, J.S., Hoffer, T. \myand Kilgore, S.B. (1982).
\textit{High School Achievements: Public, Catholic and Other Private Schools
Compared.} New York: Basic.

\bib
Cordeiro, G.M. \myand McCullagh, P. (1991). Bias correction in generalized linear 
models. \textit{Journal of the Royal Statistical Society, Series B}, 
\textbf{53}, 629--643.

\bib
Faraway, J.J. (2016). \textit{Extending the Linear Model with \textsf{R}. Second Edition.}
Boca Raton, Florida: CRC Press.

\bib
Gradshteyn, I.S. \myand Ryzhik, I.M. (1994).
\textit{Tables of Integrals, Series, and Products, Fifth Edition.}
San Diego, California: Academic Press.

\bib
Jiang, J. \myand Nguyen, T. (2017). \textit{Linear and Generalized Linear Mixed Models
and Their Applications, Second Edition.} New York: Springer.

\bib
Jiang, J., Wand, M.P. \myand Bhaskaran, A. (2022).
Usable and precise asymptotics for generalized linear
mixed model analysis and design.
\textit{Journal of the Royal Statistical Society, Series B}, 
\textbf{84}, 55--82.

\bib
Jo, S. \myand Lee, W. (2017).
Comparing the efficiency of dispersion parameter estimators
in gamma generalized linear models. \textit{The Korean Journal of Applied Statistics}, 
\textbf{30}, 95--102.

\bib
J{\o}rgensen, B. (1987). Exponential dispersion models (with discussion). 
\textit{Journal of the Royal Statistical Society, Series B}, \textbf{49},
127--162.

\bib
McCullagh, P. and Nelder, J.A. (1989). {\it Generalized Linear
Models, Second Edition}. London: Chapman and Hall.

\bib
McCulloch, C.E., Searle, S.R. \myand Neuhaus, J.M. (2008). 
\textit{Generalized, Linear, and Mixed Models. Second Edition.}
New York: John Wiley \& Sons.

\bib
Nelder, J.A. \myand Mead, R. (1965). A simplex method for function
minimization. \textit{The Computational Journal}, \textbf{7}, 308--313.

\bib
Pinheiro, J., Bates, D., \textsf{R} Core Team (2022). 
\textsf{nlme}: linear and nonlinear mixed effects models. 
\textsf{R} package version 3.1-157.
\texttt{https://CRAN.R-project.org/package=nlme}

\bib
\textsf{R} Core Team (2022). \textsf{R}: A language and environment for
statistical computing. \textsf{R} Foundation for Statistical Computing,
Vienna, Austria. \texttt{https://www.R-project.org}

\bib
Stroup, W.W. (2013). \textit{Generalized Linear Mixed Models.}
Boca Raton, Florida: CRC Press.

\vfill\eject
%
%
\renewcommand{\theequation}{S.\arabic{equation}}
\renewcommand{\thesection}{S.\arabic{section}}
\renewcommand{\thetable}{S.\arabic{table}}
\setcounter{equation}{0}
\setcounter{table}{0}
\setcounter{section}{0}
\setcounter{page}{1}
\setcounter{footnote}{0}

\begin{center}

{\Large Supplement for:}
\vskip3mm

\centerline{\Large\bf Dispersion Parameter Extension of Generalized}
\vskip1mm
\centerline{\Large\bf Linear Mixed Model Asymptotics}
\vskip7mm
\ifthenelse{\boolean{UnBlinded}}{\centerline{\normalsize\sc Aishwarya Bhaskaran and Matt P. Wand}
\vskip5mm
\centerline{\textit{University of Technology Sydney}}
\vskip6mm}{\null}
\end{center}

\section{Introduction}

This supplement contains derivational details concerning the article's results.
In Section \ref{sec:theProof} we provide a proof of Theorem \ref{thm:mainTheorem}.
Section \ref{sec:IllusSimplic} provides illustration of how, for a special case, the (A3) moment 
condition can be reduced to a simpler moment condition.

\section{Proof of Theorem \ref{thm:mainTheorem}}\label{sec:theProof}

In this section we provide a proof of Theorem \ref{thm:mainTheorem}. 
We start by setting up notation for some key quantities that arise throughout
the proof, as well as useful generic mathematical notation. The main body 
of the proof involves asymptotic approximation of the Fisher information
matrix and its inverse. Much of this was achieved in Appendix A of Jiang \textit{et al.} (2022).
However, the dispersion parameter extension leads to new Fisher information entries.
We then apply Lemma 2 of Jiang \textit{et al.} (2022) to establish an asymptotic equivalence
result between the matrix square roots of two relevant approximations
to the inverse Fisher information matrix. The final steps required to establish 
Theorem \ref{thm:mainTheorem} are then carried out.

\subsection{Notation}

We divide the notation into two parts: (1) that for key quantities specific to the
model at hand and (2) generic mathematical notation that aids the proof.

\subsubsection{Notation for Some Key Quantities}

Throughout this proof we let 
$$\psi\equiv1/\phi$$ 
denote the \emph{reciprocal} dispersion parameter. Working with $\psi$, rather than
$\phi$, in Fisher information approximations involves simpler expressions in the 
derivation of the asymptotic joint normality result for the model parameters.
The transformation from $\psi$ to $\phi=1/\psi$, using the Multivariate Delta Method,
is carried out after such a result is established.

For each $1\le i\le m$ and $1\le j\le n_i$ let $\bX_{ij}\equiv(\bXAij^T,\bXBij^T)^T$ and
$\bX_i\equiv(\bX_{i1},\ldots,\bX_{in_i})$. Let $\bY_i$, $1\le i\le m$, be defined 
analogously.

Define $\GscrAi$ and $\HscrAAi$, for each $1\le i\le m$, as follows
{\setlength\arraycolsep{1pt}
\begin{eqnarray*}
\GscrAi&\equiv&\dispsumjni\{Y_{ij} -b'\big(\linXABUi\big)\}\bXAij,\\[1ex]
\HscrAAi&\equiv&\dispsumjni b''\big(\linXABUi\big)\bXAij\bXAij^T,
\end{eqnarray*}
}
In a similar vein, define $\HscrdAAAi$ to be the $\dA\times\dA\times\dA$ 
array with $(r,s,t)$ entry equal to 
$$\dispsumjn b'''\big(\linXABUi\big)(\bXAij)_r(\bXAij)_s(\bXAij)_t$$

\subsubsection{Generic Mathematical Notation}

For a generic $d\times1$ vector $\bv$ we define $\bv^{\otimes 2}\equiv \bv\bv^T$.
We also let $\diag(\bv)$ denote the $d\times d$ diagonal matrix
with the entries of $\bv$ along the diagonal. For a matrix $\bM$
let $\Vert\bM\VertF=\{\tr(\bM^T\bM)\}^{1/2}$ denote its Frobenius norm.

For $f$ a smooth real-valued function of the $d$-variate argument 
$\bx\equiv(x_1,\ldots,x_d)$, let $\nabla f(\bx)$ denote
the $d\times1$ vector with $i$th entry $\partial f(\bx)/\partial x_i$,
$\nabla^2 f(\bx)$ denote the $d\times d$ matrix with $(i,j)$ entry 
$\partial^2 f(\bx)/(\partial x_i\partial x_j)$ and 
$\nabla^3 f(\bx)$ denote the $d\times d\times d$ array with $(i,j,k)$ entry 
$\partial^3 f(\bx)/(\partial x_i\partial x_j\partial x_k)$.

If $\bAsc$ is a $d_1\times d_2\times d_3$ array and $\bM$ is a $d_1\times d_2$
matrix then we let
$$\bAsc\bigstar\bM\quad\mbox{denote the $d_3\times1$ vector with $t$th entry given by}
\quad\sum_{r=1}^{d_1}\sum_{s=1}^{d_2}(\bAsc)_{rst}\bM_{rs}.
$$

\subsection{Fisher Information Approximation}

The Fisher information corresponding to the parameter vector 
\begin{equation}
(\bbetaA,\bbetaB,\vech(\bSigma),\psi)
\label{eq:FirstCapitol}
\end{equation}
is a symmetric matrix having $\dA+\dB+\smhalf\dA(\dA+1)+1$ rows and columns.
Appendix A of Jiang \textit{et al.} (2022) provides an adequate approximation of the
$\big(\bbetaA,\bbetaB,\vech(\bSigma)\big)$ diagonal block. The dispersion
parameter extension requires similar approximations for the $\psi$ diagonal
block and the $\big(\bbetaA,\bbetaB,\vech(\bSigma),\psi\big)$ off-diagonal block.
The sub-blocks of the $\big(\bbetaA,\bbetaB,\vech(\bSigma),\psi\big)$ off-diagonal block
correspond to each of 
$$\quad\big(\bbetaA,\bbetaB,\vech(\bSigma),\big),
\quad\big(\bbetaA,\bbetaB,\psi\big)
\quad\mbox{and}\quad\big(\vech(\bSigma),\psi\big)
$$
The first of these is treated in Jiang \textit{et al.} (2022).
The second and third of these are treated in Sections \ref{sec:floaties} and \ref{sec:SunnyBoy}
respectively.

\subsubsection{Higher Order Approximation of Multivariate Integral Ratios}\label{sec:Miyata}

Our main tool for approximation of the Fisher information matrix entries
for generalized linear mixed models is higher order Laplace-type approximation of 
multivariate integral ratios. Appendix A of Miyata (2004) provides
such a result, which states that for smooth real-valued $d$-variate 
functions $\gothicg$, $\gothicc$ and $\gothich$ we have
\begin{equation}
{\setlength\arraycolsep{1pt}
\begin{array}{rcl}
&&{\displaystyle\frac{\int_{\real^d}\gothicg(\bx)\gothicc(\bx)\exp\{-n \gothich(\bx)\}\,d\bx}
{\int_{\real^d} \gothicc(\bx)\exp\{-n \gothich(\bx)\}\,d\bx}}
=\gothicg(\bx^*)+
{\displaystyle\frac{\nabla \gothicg(\bx^*)^T\{\nabla^2 \gothich(\bx^*)\}^{-1}\nabla \gothicc(\bx^*)}{n \gothicc(\bx^*)}}\\[3ex]
&&\qquad+{\displaystyle\frac{\tr[\{\nabla^2 \gothich(\bx^*)\}^{-1}\nabla^2 \gothicg(\bx^*)]}{2n}}
-{\displaystyle\frac{\nabla \gothicg(\bx^*)^T\Big[\nabla^3 \gothich(\bx^*)\bigstar\{\nabla^2 \gothich(\bx^*)\}^{-1}\Big]}
{2n}}+O(n^{-2})
\end{array}
}
\label{eq:Miyata}
\end{equation}
where 
$$\bx^*\equiv\argmin{\bx\in\real^d}\gothich(\bx).$$

\subsubsection{The $\bU_i^*$ Quantity and Its Approximation}

For all required Fisher information approximations for the $i$th group, 
the $h$ function appearing in (\ref{eq:Miyata}) corresponds to the stochastic 
function
\begin{equation}
\gothich_i(\bu)\equiv\,-\frac{\psi}{n}\dispsumjni\left\{Y_{ij}\bu^T\bXAij-b\big(\linXABu\big)\right\}
\label{eq:BuzzAldrin}
\end{equation}
and its minimum is denoted by the random vector
$$\bU_i^*\equiv\argmin{\bu\in\real^d}\gothich_i(\bu).$$
Taylor series expansion, similar to that given in Appendix A.3.1 of Jiang \textit{et al.} (2022),
followed by asymptotic series inversion leads to the three-term approximation:
\begin{equation}
\bU_i^*=\bU_i+\HscrAAi^{-1}\GscrAi-\smhalf\HscrAAi^{-1}
\Big\{\HscrdAAAi\bigstar\big(\HscrAAi^{-1}\GscrAi\GscrAi^T\HscrAAi^{-1}\big)\Big\}
+O_P(n^{-3/2})\bone_{\dA}.
\label{eq:UiStar}
\end{equation}

\subsubsection{The $\big(\bbetaA,\bbetaB,\vech(\bSigma)\big)$ Diagonal Block}

For the case of the dispersion parameter  being fixed rather than estimated, 
Jiang \textit{et al.} (2022) derives an approximation to the Fisher information of $
\big(\bbetaA,\bbetaB,\vech(\bSigma)\big)$. Appendix A.5 provides the resultant
approximation. For the extension to $\phi=1/\psi$ being estimated, this approximation
corresponds to the $\big(\bbetaA,\bbetaB,\vech(\bSigma)\big)$ diagonal block of
the Fisher information matrix for the extended parameter vector (\ref{eq:FirstCapitol}).

\subsubsection{The $\psi$ Diagonal Block}

The $i$th contribution to the score is
{\setlength\arraycolsep{1pt}
\begin{eqnarray*}
\frac{\partial\log\{\pYiGXi(\bY_i|\bX_i)\}}{\partial\psi}&=&
\dispsumjni\left[\Big\{Y_{ij}\big(\linXAB\big)+\gothicc(Y_{ij})\Big\}
-\frac{d d\big(1/\psi\big)}{d\psi}\right]\\
&&\qquad+\frac{{\displaystyle\intdA} \gothicg_i(\bu)\gothicc(\bu)\exp\{\psi \gothicg_i(\bu)\}\,d\bu}
{{\displaystyle\intdA}\gothicc(\bu)\exp\{\psi \gothicg_i(\bu)\}\,d\bu}.
\end{eqnarray*}
}
The derivative of the $i$th contribution to the score is
\begin{equation}
\frac{\partial^2\log\{\pYiGXi(\bY_i|\bX_i)\}}{\partial\psi^2}=
-\frac{n_i\,d^2 d\big(1/\psi\big)}{d\psi^2}+\Qsc_{1i}
-\Qsc_{2i}^2
\label{eq:DontLookBack}
\end{equation}
where
\begin{equation}
\Qsc_{1i}\equiv 
\frac{{\displaystyle\intdA} \gothicg_i^2(\bu)\gothicc(\bu)\exp\{-n \gothich_i(\bu)\}\,d\bu}
{{\displaystyle\intdA} \gothicc(\bu)\exp\{-n \gothich_i(\bu)\}\,d\bu}
\label{eq:AshleyMallett}
\end{equation}
and
\begin{equation}
\Qsc_{2i}\equiv 
\frac{{\displaystyle\intdA} \gothicg_i(\bu)\gothicc(\bu)\exp\{-n \gothich_i(\bu)\}\,d\bu}
{{\displaystyle\intdA} \gothicc(\bu)\exp\{-n \gothich_i(\bu)\}\,d\bu}
\label{eq:TerryJenner}
\end{equation}
with
$$\gothicc(\bu)\equiv\exp\left(-\smhalf\bu^T\bSigma^{-1}\bu\right),\quad
\gothicg_i(\bu)\equiv\dispsumjni\left\{Y_{ij}\bu^T\bXAij-b\big(\linXABu\big)\right\}
$$
and $\gothich_i$ is as given by (\ref{eq:BuzzAldrin}). Application of (\ref{eq:Miyata}) to each
of (\ref{eq:AshleyMallett}) and (\ref{eq:TerryJenner}) and use of (\ref{eq:UiStar}) 
leads to the following three-term approximations to $\Qsc_{1i}$ and $\Qsc_{2i}$:
$$\Qsc_{1i}=\gothicg_i(\bU_i)^2+\GscrAi^T\HscrAAi^{-1}\GscrAi\,\gothicg_i(\bU_i)
-\frac{\dA\,\gothicg_i(\bU_i)}{\psi}+O_P(n^{1/2})
$$
and
$$\Qsc_{2i}=\gothicg_i(\bU_i)+\smhalf\GscrAi^T\HscrAAi^{-1}\GscrAi-\frac{\dA}{2\psi}+O_P(n^{-1/2}).$$
Therefore
{\setlength\arraycolsep{1pt}
\begin{eqnarray*}
\Qsc_{1i}-\Qsc_{2i}^2&=&\gothicg_i(\bU_i)^2+\GscrAi^T\HscrAAi^{-1}\GscrAi\,\gothicg_i(\bU_i)
-\frac{\dA \gothicg_i(\bU_i)}{\psi}+O_P(n^{1/2})\\
&&\quad-\left\{\gothicg_i(\bU_i)+\smhalf\GscrAi^T\HscrAAi^{-1}\GscrAi-\frac{\dA}{2\psi}+O_P(n^{-1/2})\right\}^2\\[1ex]
&=&\gothicg_i(\bU_i)^2+\GscrAi^T\HscrAAi^{-1}\GscrAi\,\gothicg_i(\bU_i)
-\frac{\dA \gothicg_i(\bU_i)}{\psi}+O_P(n^{1/2})\\
&&\quad-\left\{\gothicg_i(\bU_i)^2+\GscrAi^T\HscrAAi^{-1}\GscrAi\,\gothicg_i(\bU_i)
-\frac{\dA \gothicg_i(\bU_i)}{\psi}\right\}+O_P(n^{1/2})\\[1ex]
&=&O_P(n^{1/2})
\end{eqnarray*}
}
and we arrive at the approximation
$$
\frac{\partial^2\log\{\pYiGXi(\bY_i|\bX_i)\}}{\partial\psi^2}=
-\frac{n_i\,d^2 d\big(1/\psi\big)}{d\psi^2}+O_P(n^{1/2}).
$$
Therefore, the $\psi$ diagonal block of the Fisher information is 
$$mn\,\left(\frac{d^2 d\big(1/\psi\big)}{d\psi^2}\right)+O_P(mn^{1/2}).
$$

\subsubsection{The $\big(\bbetaA,\bbetaB,\psi\big)$ Off-Diagonal Block}\label{sec:floaties}

Let 
$$\bbeta\equiv\left[
\begin{array}{c}
\bbetaA\\[1ex]
\bbetaB
\end{array}
\right]
$$
denote the full fixed effects vector. As established in the appendix of Wand (2007), 
the $i$th contribution to the partial derivative, with respect to $\psi$, of the $\bbeta$ score is
$$\frac{\partial\delbbeta\log\{\pYiGXi(\bY_i|\bX_i)\}}{\partial\psi}=\bX_i^T[\bY_i-E\{E(\bY_i|\bU_i)|\bY_i\}].$$
Noting that
{\setlength\arraycolsep{1pt}
\begin{eqnarray*}
E\big[\bX_i^T\{\bY_i-E\{E(\bY_i|\bU_i)|\bY_i\}|\bX_i\big]
&=&\bX_i^T\Big(E(\bY_i)-E[E\{E(\bY_i|\bU_i)|\bY_i\}]\Big)\\[1ex]
&=&\bX_i^T\big\{E(\bY_i)-E(\bY_i)\big\}=\bO
\end{eqnarray*}
}
it is apparent that the $\big(\bbetaA,\bbetaB,\psi\big)$ off-diagonal block of the Fisher information
has all entries being exactly zero.

\subsubsection{The $\big(\vech(\bSigma),\psi\big)$ Off-Diagonal Block}\label{sec:SunnyBoy}

The $i$th contribution to the second order partial derivative with respect to $\vech(\bSigma)$ and $\psi$ is
\begin{equation}
{\setlength\arraycolsep{1pt}
\begin{array}{rcl}
\delVechSigma{\displaystyle\frac{\partial\log\{\pYiGXi(\bY_i|\bX_i)\}}{\partial\psi}}
&=&\delVechSigma\left[{\displaystyle\frac{{\displaystyle\intdA} \gothicg_i(\bu)
\gothicc(\bu)\exp\{-n \gothich_i(\bu)\}\,d\bu}
{{\displaystyle\intdA}\gothicc(\bu)\exp\{-n \gothich_i(\bu)\}\,d\bu}}\right]\\[5ex]
&=&\smhalf\Big(\bQsc_{3i}-\Qsc_{2i}\bQsc_{4i}\Big)
\end{array}
}
\label{eq:CarefulWithAxe}
\end{equation}
where
\begin{equation}
\bQsc_{3i}\equiv 
\frac{{\displaystyle\intdA} \gothicg_i(\bu)
\bD_{\dA}^T\vecof(\bSigma^{-1}\bu\bu^T\bSigma^{-1})\gothicc(\bu)\exp\{-n \gothich_i(\bu)\}\,d\bu}
{{\displaystyle\intdA} \gothicc(\bu)\exp\{-n \gothich_i(\bu)\}\,d\bu}
\label{eq:ShaneWarne}
\end{equation}
and
\begin{equation}
\bQsc_{4i}\equiv 
\frac{{\displaystyle\intdA} 
\bD_{\dA}^T\vecof(\bSigma^{-1}\bu\bu^T\bSigma^{-1})\gothicc(\bu)\exp\{-n \gothich_i(\bu)\}\,d\bu}
{{\displaystyle\intdA} \gothicc(\bu)\exp\{-n \gothich_i(\bu)\}\,d\bu}.
\label{eq:StuartMacGill}
\end{equation}
Application of (\ref{eq:Miyata}) to each of (\ref{eq:ShaneWarne}) 
and (\ref{eq:StuartMacGill}) and use of (\ref{eq:UiStar}) 
leads to the following approximations to $\bQsc_{3i}$ and $\bQsc_{4i}$:
$$\bQsc_{3i}=\bD_{\dA}^T\vecof\Big(\bSigma^{-1}(\bU_i\bU_i^T + 2\HscrAAi^{-1}\GscrAi\bU_i^T)\bSigma^{-1}\Big)
\gothicg_i(\bU_i)+O_P(1)\bone_{\dAwindow}$$
and
$$\bQsc_{4i}=\bD_{\dA}^T\vecof\Big(\bSigma^{-1}(\bU_i\bU_i^T + 2\HscrAAi^{-1}\GscrAi\bU_i^T)\bSigma^{-1}\Big)
+O_P(n^{-1})\bone_{\dAwindow}
$$
where $\dAwindow\equiv\dA(\dA+1)/2$. 
Substitution of these approximations into (\ref{eq:CarefulWithAxe}) then gives
{\setlength\arraycolsep{1pt}
\begin{eqnarray*}
\delVechSigma\frac{\partial\log\{\pYiGXi(\bY_i|\bX_i)\}}{\partial\psi}
&=&\smhalf\bD_{\dA}^T\vecof\Big(\bSigma^{-1}(\bU_i\bU_i^T + 2\HscrAAi^{-1}\GscrAi\bU_i^T)\bSigma^{-1}\Big)
\gothicg_i(\bU_i)\\[1ex]
&&\quad -\smhalf\Big(\gothicg_i(\bU_i)+O_P(1)\Big)\Big\{\bD_{\dA}^T\vecof\Big(\bSigma^{-1}(\bU_i\bU_i^T\\
&&\quad\qquad+ 2\HscrAAi^{-1}\GscrAi\bU_i^T)\bSigma^{-1}\Big)
+O_P(n^{-1})\bone_{\dAwindow}\Big\}\\[1ex]
&=&O_P(1)\bone_{\dAwindow}.
\end{eqnarray*}
}
It follows that the $\big(\vech(\bSigma),\psi\big)$ off-diagonal block of the Fisher
information matrix is
$$-E\left[\sumim\delVechSigma\frac{\partial\log\{\pYiGXi(\bY_i|\bX_i)\}}{\partial\psi}\right]=O_P(m)\bone_{\dAwindow}.$$

\subsubsection{Assembly of Fisher Information Sub-Block Approximations}

From the Fisher information sub-block approximations obtained in the previous 
six subsubsections we have
{\setlength\arraycolsep{1pt}
\begin{eqnarray*}
&&I\Big(\bbetaA,\bbetaB,\vech(\bSigma),\psi\Big)\\[1ex]
&&=
\left[
{\setlength\arraycolsep{0pt}
\begin{array}{cccc}
m\bSigma^{-1}+O_P(mn^{-1})\bone_{\dA}^{\otimes 2} 
& O_P(m)\bone_{\dA}\bone_{\dB}^T & O_P(mn^{-1})\bone_{\dA}\bone_{\dAwindow}^T & \bzero_{\dA}\\[3ex]
O_P(m)\bone_{\dB}\bone_{\dA}^T
&\displaystyle{\frac{mn\bLambda_{\bbetaB}^{-1}}{\phi}}+o_P(mn)\bone_{\dB}^{\otimes2}
&O_P(m)\bone_{\dB}\bone_{\dAwindow}^T & \bzero_{\dB}\\[3ex]
O_P(mn^{-1})\bone_{\dAwindow}\bone_{\dA}^T & O_P(m)\bone_{\dAwindow}\bone_{\dB}^T  
&\displaystyle{\frac{m\bD_{\dA}^T(\bSigma^{-1}\otimes\bSigma^{-1})\bD_{\dA}}{2}} 
& O_P(m)\bone_{\dAwindow} \\
&&+O_P(mn^{-1})\bone_{\dAwindow}^{\otimes2}\\[3ex]
\bzero_{\dA}^T &  \bzero_{\dB}^T    & O_P(m)\bone_{\dAwindow}^T   
& mn\left(\displaystyle{\frac{d^2d(1/\psi)}{d\psi^2}}\right)\\[2.5ex]
          &         &      & + O_P(mn^{1/2})
\end{array}
}
\right].
\end{eqnarray*}
}

\subsection{Inverse Fisher Information Approximation}

Note that
$$I\Big(\bbetaA,\bbetaB,\vech(\bSigma),\psi\Big)
=\left[
\begin{array}{cc}
\bA_{11}   &  \bA_{12} \\[1ex]
\bA_{12}^T &  A_{22}
\end{array}
\right]
$$
where
\begin{equation}
\bA_{11}\equiv
\left[
{\setlength\arraycolsep{2pt}
\begin{array}{ccc}
m\bSigma^{-1}+O_P(mn^{-1})\bone_{\dA}^{\otimes 2} 
& O_P(m)\bone_{\dA}\bone_{\dB}^T & O_P(mn^{-1})\bone_{\dA}\bone_{\dAwindow}^T\\[3ex]
O_P(m)\bone_{\dB}\bone_{\dA}^T
&\displaystyle{\frac{mn\bLambda_{\bbetaB}^{-1}}{\phi}}+o_P(mn)\bone_{\dB}^{\otimes2}
&O_P(m)\bone_{\dB}\bone_{\dAwindow}^T\\[3ex]
O_P(mn^{-1})\bone_{\dAwindow}\bone_{\dA}^T & O_P(m)\bone_{\dAwindow}\bone_{\dB}^T  
&\displaystyle{\frac{m\bD_{\dA}^T(\bSigma^{-1}\otimes\bSigma^{-1})\bD_{\dA}}{2}}\\
&&+O_P(mn^{-1})\bone_{\dAwindow}^{\otimes2}
\end{array}
}
\right],
\label{eq:UCIrvine}
\end{equation}
$$A_{22}\equiv mn E\left\{\displaystyle{\frac{d^2d(\frac{1}{\psi})}{d\psi^2}}\right\}
+O_P(mn^{1/2})
\quad\mbox{and}\quad
\bA_{12}\equiv\big[\bzero_{\dA}^T \ \ \ \bzero_{\dB}^T\ \ \ O_P(m)\bone_{\dAwindow}^T\big]^T.$$
Let
$$I\Big(\bbetaA,\bbetaB,\vech(\bSigma),\psi\Big)^{-1}
=\left[
\begin{array}{cc}
\bA^{11}   &  \bA^{12} \\[1ex]
(\bA^{12})^T &  A^{22}
\end{array}
\right].
$$
The upper left block of $I\Big(\bbetaA,\bbetaB,\vech(\bSigma),\psi\Big)^{-1}$ is
\begin{equation}
\bA^{11}=\bA_{11}^{-1}+\frac{\bA_{11}^{-1}\bA_{12}\bA_{12}^T\bA_{11}^{-1}}
{A_{22}-\bA_{12}^T\bA_{11}^{-1}\bA_{12}}.
\label{eq:CockOfTheNorth}
\end{equation}
The leading terms of $\bA_{11}^{-1}$ are provided by equation (A27) of Jiang \textit{et al.} (2022).
Since $\bA_{12}^T\bA_{11}^{-1}\bA_{12}=O_P(m)$ we have 
$$A_{22}-\bA_{12}^T\bA_{11}^{-1}\bA_{12}\equiv 
mn\left(\displaystyle{\frac{d^2d(1/\psi)}
{d\psi^2}}\right)+O_P(mn^{1/2})
$$
and so
\begin{equation}
\frac{1}{A_{22}-\bA_{12}^T\bA_{11}^{-1}\bA_{12}}=O_P(m^{-1}n^{-1}).
\label{eq:FairyMeadow}
\end{equation}
Also,
$$\mbox{the lower right block of}\ \bA_{11}^{-1}\bA_{12}\bA_{12}^T\bA_{11}^{-1}=O_P(1)\bone_{\dAwindow}^{\otimes2}$$
and all other entries of $\bA_{11}^{-1}\bA_{12}\bA_{12}^T\bA_{11}^{-1}$ are exactly zero. Hence,
in view of (\ref{eq:FairyMeadow}), the second term of (\ref{eq:CockOfTheNorth}) is a matrix with
all entries either $O_P(m^{-1}n^{-1})$ or zero. Consequently, $\bA^{11}$ equals the right-hand
side of (A27) of Jiang \textit{et al.} (2022) except for a rearrangement of the entries to concur
with the $\big(\bbetaA,\bbetaB,\vech(\bSigma)\big)$ parameter ordering rather than 
$\big(\bbetaA,\vech(\bSigma),\bbetaB\big)$.

The lower right entry of 
$$I\Big(\bbetaA,\bbetaB,\vech(\bSigma),\psi\Big)^{-1}$$ 
is
\begin{equation}
A^{22}=\frac{1}{A_{22}}+\frac{\bA_{12}^T
(\bA_{11}-\bA_{12}\bA_{12}^T/A_{22})^{-1}\bA_{12}}{A_{22}^2}.
\label{eq:GabbaTest}
\end{equation}
Note that 
$$\frac{1}{A_{22}}=\left(\displaystyle{
\frac{d^2d(1/\psi)}{d\psi^2}}\right)^{-1}(mn)^{-1}+O_P(m^{-1}n^{-3/2}).
$$
Also
$$\mbox{the lower right block of}\ \bA_{12}\bA_{12}^T/A_{22}=O_P(mn^{-1})\bone_{\dAwindow}^{\otimes2}$$
and all other entries of $\bA_{12}\bA_{12}^T/A_{22}$ are exactly zero. It follows that the expression
on the right-hand side of (\ref{eq:UCIrvine}) also holds for $\bA_{11}-\bA_{12}\bA_{12}^T/A_{22}$
and, using equation (A27) of Jiang \textit{et al.} (2022),
$$\frac{\bA_{12}^T(\bA_{11}-\bA_{12}\bA_{12}^T/A_{22})^{-1}\bA_{12}}{A_{22}^2}=O_P(m^{-1}n^{-2}).$$
Hence, the second term on the right-hand side of (\ref{eq:GabbaTest}) is asymptotically negligible 
compared with the first term and we have
$$A^{22}=\left(\displaystyle{\frac{d^2d(1/\psi)}{d\psi^2}}\right)^{-1}
(mn)^{-1}+O_P(m^{-1}n^{-3/2}).
$$
The off-diagonal block of the inverse Fisher information matrix is
$$\bA^{12}=-(\bA_{11}-\bA_{12}\bA_{12}^T/A_{22})^{-1}\bA_{12}/A_{22}
=\Big[\bzero_{\dA}^T \ \ \ \bzero_{\dB}^T\ \ \ O_P(m^{-1}n^{-1})\bone_{\dAwindow}^T\Big]^T$$
and so we have 
\begin{equation}
{\setlength\arraycolsep{1pt}
\begin{array}{rcl}
&&I\big(\bbetaA,\bbetaB,\vechbSigma,\psi\big)^{-1}
=I\big(\bbetaA,\bbetaB,\vechbSigma,\psi\big)_{\infty}^{-1}\\[2ex]
&&\qquad\qquad +{\displaystyle\frac{1}{mn}}
\left[
{\setlength\arraycolsep{3pt}
\begin{array}{cccc}
O_P(1)\bone_{\dA}^{\otimes2} & O_P(1)\bone_{\dA}\bone_{\dB}^T  & O_P(1)\bone_{\dA}\bone_{\dAwindow}^T 
& \bzero_{\dA} \\[1.5ex]
O_P(1)\bone_{\dB}\bone_{\dA}^T &o_P(1)\bone_{\dB}^{\otimes 2}& O_P(1) \bone_{\dB}\bone_{\dAwindow}^T  
& \bzero_{\dB} \\[1.5ex]
O_P(1)\bone_{\dAwindow}\bone_{\dA}^T & O_P(1)\bone_{\dAwindow}\bone_{\dB}^T &  O_P(1)\bone_{\dAwindow}^{\otimes2} 
& O_P(1)\bone_{\dAwindow}\\[1.5ex]
\bzero_{\dA}^T  &  \bzero_{\dB}^T &  O_P(1)\bone_{\dAwindow}^T & O_P(n^{-1/2}) 
\end{array}
}
\right]
\end{array}
}
\label{eq:actualInfo}
\end{equation}
where
$$I\big(\bbetaA,\bbetaB,\vechbSigma,\psi\big)_{\infty}^{-1}
\equiv
\left[
\begin{array}{cccc}
{\displaystyle\frac{\bSigma}{m}} & \qquad \bO\qquad &\qquad\bO\qquad  & \bO\\[1ex]
\bO &\quad {\displaystyle\frac{\phi\bLambda_{\bbetaB}}{mn}} & \bO & \bO\\[1ex]
\bO & \bO & {\displaystyle\frac{2\bD_{\dA}^{+}(\bSigma\otimes\bSigma)\bD_{\dA}^{+T}}{m}} & \bO \\[2ex]
\bO & \bO &\bO  & {\displaystyle\frac{\left(\displaystyle{
\frac{d^2d(1/\psi)}{d\psi^2}}\right)^{-1}}{mn}}
\end{array}
\right].
$$

\subsection{Asymptotic Equivalence of $\{I\big(\bbetaA,\bbetaB,\vechbSigma,\psi\big)^{-1}\}^{1/2}$\\
and $\{I\big(\bbetaA,\bbetaB,\vechbSigma,\psi\big)_{\infty}^{-1}\}^{1/2}$
}

Our aim in this subsection is to establish asymptotic equivalence between 
$\{I\big(\bbetaA,\bbetaB,\vechbSigma,\psi\big)^{-1}\}^{1/2}$
and $\{I\big(\bbetaA,\bbetaB,\vechbSigma,\psi\big)_{\infty}^{-1}\}^{1/2}$
in the following sense:
\begin{equation}
\Big\Vert\{I\big(\bbetaA,\bbetaB,\vechbSigma,\psi\big)_{\infty}^{-1}\}^{-1/2}
\{I\big(\bbetaA,\bbetaB,\vechbSigma,\psi\big)^{-1}\}^{1/2}-\bI\BigVertF\convprob 0.
\label{eq:CrocsDandJ}
\end{equation}
Without loss of generality, we change the ordering of the parameters from 
$\big(\bbetaA,\bbetaB,\vechbSigma,\psi\big)$ to $\big(\bbetaA,\vechbSigma,\bbetaB,\psi\big)$ 
and note that
$$I\big(\bbetaA,\vechbSigma,\bbetaB,\psi\big)_{\infty}^{-1}=
\frac{1}{m}\left[
\begin{array}{cc}
\bK & \bO \\[1ex]
\bO & \oon\bL
\end{array}
\right]
$$
where 
$$\bK\equiv
\left[
\begin{array}{cc}
\bSigma & \bO \\[1ex]
\bO     & 2\bD_{\dA}^{+}(\bSigma\otimes\bSigma)\bD_{\dA}^{+T}
\end{array}
\right]
\qquad\mbox{and}\qquad
\bL\equiv
\left[
\begin{array}{cc}
\phi\bLambda_{\bbetaB}  & \bO \\[1ex]
\bO     & \left(\displaystyle{
\frac{d^2d(1/\psi)}{d\psi^2}}\right)^{-1}
\end{array}
\right].
$$
Also, 
$$I\big(\bbetaA,\vechbSigma,\bbetaB,\psi\big)^{-1}=
\frac{1}{m}\left[
\begin{array}{cc}
\bK+O_p(n^{-1})\bone_{\dA+\dAwindow}^{\otimes2}   & O_p(n^{-1})\bone_{\dA+\dAwindow}\bone_{\dB+1}^T \\[2ex]
O_P(n^{-1})\bone_{\dB+1}^T\bone_{\dA+\dAwindow}   & \frac{1}{n}\bL+o_p(n^{-1})\bone_{\dB+1}^{\otimes2}
\end{array}
\right]
$$
and (\ref{eq:CrocsDandJ}) follows from Lemma 2 of Jiang \textit{et al.} (2022).

\subsection{Final Steps}

Using steps analogous to those given in Appendix A.8 of Jiang \textit{et al.} (2022) we have
$$
\sqrt{m}\left[
{\setlength\arraycolsep{0pt}
\begin{array}{c}
\bbetaAMLE-\bbetaAzero\\[1.5ex]
\sqrt{n}\Big(\bbetaBMLE-\bbetaBzero\Big)\\[1.5ex]
\vech(\bSigmaMLE-\bSigmaZero)\\[1.5ex]
\sqrt{n}(\psiMLE-\psiZero)
\end{array}
}
\right]\convdist N\left(
\left[
{\setlength\arraycolsep{0pt}
\begin{array}{c}
\bzero\\[2.5ex]
\bzero\\[2.5ex]
\bzero\\[2.5ex]
0
\end{array}
}
\right],
\left[
{\setlength\arraycolsep{0pt}
\begin{array}{cccc}
\bSigmaZero & \bO                    &    \bO  & \bO \\[1ex]
\bO         & \phiZero\bLambda_{\bbetaB} &    \bO  & \bO\\[1ex]
\bO         & \bO                    & 2\bD_{\dA}^{+}(\bSigmaZero\otimes\bSigmaZero)\bD_{\dA}^{+T} & \bO\\[1ex]
\bO         & \bO                    & \bO  & \left(\displaystyle{
\frac{d^2d(1/\psi)}{d\psi^2}}\right)^{-1}_{\psi=\psiZero}
\end{array}
}
\right]
\right).
$$
Application of the Multivariate Delta Method (e.g.\ Agresti, 2013, Section 16.1.3) with the mapping
$$\Big(x_1,\ldots,x_{\dA+\dB+\dAwindow},x_{\dA+\dB+\dAwindow+1}\Big)\mapsto
\Big(x_1,\ldots,x_{\dA+\dB+\dAwindow},1/x_{\dA+\dB+\dAwindow+1}\Big)
$$
leads to Theorem \ref{thm:mainTheorem}.

\section{Illustration of Simplification of (A3)}\label{sec:IllusSimplic}

Assumption (A3) of Theorem \ref{thm:mainTheorem} is that 
\begin{equation}
\displaystyle{E\left[\frac{
E\Big[\max\big(1,\Vert\bX\Vert\big)^8\,\max\big\{1,b''\big(\linXABU\big)\big\}^4\Big|\bU\Big]}
{\min\big\{1,\lambdaMin\big(E\{\bXA\bXA^T
\,b''\big(\linXABU\big)|\bU\}\big)\big\}^2}\right]}<\infty
\label{eq:WimbledonCommon}
\end{equation}
for all $\bbetaA\in\real^{\dA}$, $\bbetaB\in\real^{\dB}$ and $\bSigma$ a $\dA\times\dA$
symmetric and positive definite matrix. 

In this section we prove that, for the special case: 
\begin{equation}
\dA=1,\quad\bXA=1\quad\mbox{and}\quad b=\exp
\label{eq:CornFlakeFactory}
\end{equation}
corresponding to Poisson responses, assumption (A3) is implied by the moment generating
function existence condition
\begin{equation}
E\{\exp(\bt^T\bXB)\}<\infty\ \mbox{for all $\bt\in\real^{\dB}$}.
\label{eq:Beserkeley}
\end{equation}
To justify the sufficiency of (\ref{eq:Beserkeley}) first note that, for the 
(\ref{eq:CornFlakeFactory}) special case, the numerator of the random variable 
inside the outermost expectation of (\ref{eq:WimbledonCommon}) equals
$$E\Big[\max\big\{1,(1+\Vert\bXB\Vert^2)^4\big\}\,
\max\big\{1,\exp\big(\betaA+U+\bbetaB^T\bXB\big)\big\}^4\big|U\Big].$$
Then application of the Cauchy-Schwartz inequality for conditional expectations gives
{\setlength\arraycolsep{1pt}
\begin{eqnarray*}
&&E\Big[\max\big\{1,(1+\Vert\bXB\Vert^2)^4\big\}\,
\max\big\{1,\exp\big(\betaA+U+\bbetaB^T\bXB\big)\big\}^4\Big|U\Big]\\[1ex]
&&\quad\le\Big(E\Big[\max\big\{1,(1+\Vert\bXB\Vert^2)^8\big\}\Big|U\Big]\Big)^{1/2}
\Big(E\big[\max\big\{1,\exp\big(\betaA+U+\bbetaB^T\bXB\big)\big\}^8\big|U\big]\Big)^{1/2}\\[1ex]
&&\quad\le\Big[1+E\big\{\big(1+\Vert\bXB\Vert^2\big)^8\big\}\Big]^{1/2}
\Big(1+E\big[\exp\big\{8(\betaA+U+\bbetaB^T\bXB)\big\}\big|U\big]\Big)^{1/2}\\[1ex]
&&\quad=\Big[1+E\big\{\big(1+\Vert\bXB\Vert^2\big)^8\big\}\Big]^{1/2}
\Big[1+\exp(8\betaA)E\big\{\exp(8\bbetaB^T\bXB)\big\}\exp(8U)\Big]^{1/2}.
\end{eqnarray*}
}
The denominator of the random variable inside the outermost expectation of (\ref{eq:WimbledonCommon}) is
{\setlength\arraycolsep{1pt}
\begin{eqnarray*}
&&\min\big\{1,E\{\exp\big(\betaA+U+\bbetaB^T\bXB\big)|U\}\big\}^2\\[1ex]
&&\qquad\qquad=\min\big\{1,\exp(2\betaA)[E\{\exp(\bbetaB^T\bXB)\}]^2\exp(2U)\}.
\end{eqnarray*}
}
Noting that, for all $x\in\real$ and $a,b>0$, 
$$\frac{\{1+a\exp(8x)\}^{1/2}}{\min\{1,b\exp(2x)\}}\le 1 
+\frac{a^{1/2}\exp(2x)}{b}+a^{1/2}\exp(4x)+\frac{\exp(-2x)}{b}$$
the random variable inside the outermost expectation of (\ref{eq:WimbledonCommon}) is bounded above by
\begin{equation}
{\setlength\arraycolsep{1pt}
\begin{array}{rcl}
&&\Big[1+E\big\{\big(1+\Vert\bXB\Vert^2\big)^8\big\}\Big]^{1/2}
\Bigg(1+\displaystyle{\frac{\exp(2\betaA)[E\{\exp(8\bbetaB^T\bXB)\}]^{1/2}\exp(2U)}{[E\{\exp(\bbetaB^T\bXB)\}]^2}}\\[2ex]
&&\qquad\qquad\qquad\qquad\qquad\qquad\qquad
+\exp(4\betaA)[E\{\exp(8\bbetaB^T\bXB)\}]^{1/2}\exp(4U)\\[1ex]
&&\qquad\qquad\qquad\qquad\qquad\qquad\qquad
+\displaystyle{\frac{\exp(-2U)}{\exp(2\betaA)[E\{\exp(\bbetaB^T\bXB)\}]^2}}\Bigg).
\end{array}
}
\label{eq:NewTown}
\end{equation}
Since $U\sim N(0,\sigsqZero)$ we have $E\{\exp(tU)\}=\exp\{\smhalf t^2\sigsqZero\}$
for all $t\in\real$. Hence, under assumption (\ref{eq:Beserkeley}), the expectation
of (\ref{eq:NewTown}) is finite which implies that (\ref{eq:WimbledonCommon}) holds
for the (\ref{eq:CornFlakeFactory}) special case. 

\section*{References}

\bib
Agresti, A. (2013). \textit{Categorical Data Analysis}. Hoboken, New Jersey: John Wiley \& Sons.

\bib
Miyata, Y. (2004).
Fully exponential Laplace approximation using asymptotic modes.
\textit{Journal of the American Statistical Association}, {\bf 99}, 
1037--1049.

\bib
Wand, M.P. (2007).
Fisher information for generalised linear mixed models.
{\it Journal of Multivariate Analysis}, {\bf 98}, 1412--1416.

%
%
%
%

\end{document}